\documentclass[12pt]{amsart}
\usepackage{amssymb}
\usepackage[all]{xypic}

\newcounter{alphalistcntr}

\theoremstyle{plain}
\newtheorem{theorem}{Theorem}[section]
\newtheorem{lemma}[theorem]{Lemma}
\newtheorem{corollary}[theorem]{Corollary}
\newtheorem{prop}[theorem]{Proposition}

\newtheorem{algorithm}[theorem]{Algorithm}

\theoremstyle{remark}

\newtheorem{remark}[theorem]{Remark}

\newtheorem*{note*}{Note}
\newtheorem*{remark*}{Remark}
\newtheorem*{example*}{Example}

\theoremstyle{definition}
\newtheorem*{definition*}{Definition}
\newtheorem{definition}[theorem]{Definition}

\newcommand{\Z}{\mathbb{Z}}

\newcommand{\Q}{\mathbb{Q}}
\newcommand{\Qu}{\mathbb{Q}}

\newcommand{\N}{\mathbb{N}}

\newcommand{\tensor}{\otimes}

\newcommand{\Cl}{\mathrm{Cl}}

\newcommand{\Aut}{\mathrm{Aut}}
\newcommand{\End}{\mathrm{End}}
\newcommand{\Mat}{\mathrm{Mat}}
\newcommand{\GL}{\mathrm{GL}}
\newcommand{\SL}{\mathrm{SL}}

\newcommand{\E}{\mathrm{E}}
\newcommand{\op}{\mathrm{op}}
\newcommand{\calA}{\mathcal{A}}
\newcommand{\calO}{\mathcal{O}}
\newcommand{\calM}{\mathcal{M}}
\newcommand{\fra}{\mathfrak{a}}
\newcommand{\frb}{\mathfrak{b}}
\newcommand{\frc}{\mathfrak{c}}

\newcommand{\frg}{\mathfrak{g}}
\newcommand{\frp}{\mathfrak{p}}
\newcommand{\ra}{\rightarrow}
\newcommand{\lra}{\longrightarrow}

\newcommand{\OL}{{\calO_L}}

\renewcommand{\OE}{{\calO_E}}
\newcommand{\ZG}{{\Z[G]}}
\newcommand{\QG}{{\Q[G]}}
\newcommand{\OEG}{{\calO_E[G]}}

\newcommand{\EG}{E[G]}

\newcommand{\sseq}{\subseteq}

\renewcommand{\labelenumi}{(\alph{enumi})}

\title[Computing generators of free modules over orders in group algebras]{Computing generators of free modules \\ over orders in group algebras}

\author{Werner Bley}
\address{Werner Bley\\
Fachbereich f\"ur Mathematik und Informatik der Universit\"at Kassel\\
Heinrich-Plett-Str. 40\\
34132 Kassel\\
Germany}
\email{bley@mathematik.uni-kassel.de}
\urladdr{http://www.mathematik.uni-kassel.de/$\sim$bley}

\author{Henri Johnston}
\address{Henri Johnston\\ 
St. Hugh's College\\
St. Margaret's Road\\
Oxford OX2 6LE\\
UK
}
\email{henri@maths.ox.ac.uk}
\urladdr{http://www.maths.ox.ac.uk/$\sim$henri}

\thanks{Johnston was supported by a grant from the
 Deutscher Akademischer Austausch Dienst.}
\date{27th January 2008}

\begin{document}

\begin{abstract}
Let $E$ be a number field and $G$ be a finite group.
Let $\mathcal{A}$ be any $\mathcal{O}_{E}$-order of full rank
in the group algebra $E[G]$ and $X$ be a (left) $\mathcal{A}$-lattice.
We give a necessary and sufficient condition for $X$ to be free of given rank $d$
over $\mathcal{A}$. In the case that the Wedderburn decomposition $E[G] \cong \oplus_{\chi} M_{\chi}$ 
is explicitly computable and each $M_{\chi}$ is in fact a matrix ring over a field,
this leads to an algorithm that either gives elements 
$\alpha_{1}, \ldots, \alpha_{d} \in X$ such that 
$X=\mathcal{A}\alpha_{1} \oplus \ldots \oplus \mathcal{A}\alpha_{d}$ 
or determines that no such elements exist.

Let $L/K$ be a finite Galois extension of number fields with Galois group $G$
such that $E$ is a subfield of $K$ and put $d=[K:E]$. 
The algorithm can be applied to certain Galois modules that arise naturally 
in this situation. For example, one can take $X$ to be $\mathcal{O}_{L}$,
the ring of algebraic integers of $L$, and $\mathcal{A}$ to be the 
associated order $\mathcal{A}(E[G];\mathcal{O}_{L}) \sseq E[G]$.
The application of the algorithm to this special situation 
is implemented in Magma under certain extra hypotheses when $K=E=\Q$.
\end{abstract}

\maketitle

\section{Introduction}\label{intro}

Let $E$ be a number field and $G$ be a finite group.
Let $\mathcal{A}$ be any $\mathcal{O}_{E}$-order of full rank in the
group algebra $E[G]$ and $X$ be a (left) $\mathcal{A}$-lattice, i.e.,
a (left) $\mathcal{A}$-module that is finitely generated and torsion-free 
over $\mathcal{O}_{E}$.
The first result of this paper is a necessary and sufficient condition for 
$X$ to be free of given rank $d$ over $\mathcal{A}$. 
In order to use this criterion
for computational purposes, we have to impose two hypotheses:

\smallskip

\begin{itemize}
\item [(H1)] The Wedderburn decomposition $E[G]\cong \oplus_{\chi} M_{\chi}$,
where each $M_{\chi}$ is a matrix ring over a division ring, is explicitly computable.
\item [(H2)] The Schur indices of all $E$-rational irreducible characters of $G$ are equal to 1,
i.e., each $M_{\chi}$ above is in fact a matrix ring over a number field.
\end{itemize}
\smallskip
Under these hypotheses, we give an algorithm that either computes elements 
$\alpha_{1}, \ldots, \alpha_{d} \in X$ such that 
$X=\mathcal{A}\alpha_{1} \oplus \ldots \oplus \mathcal{A}\alpha_{d}$ 
or determines that no such elements exist. 
More generally, 
the group algebra $E[G]$ can be replaced by any finite 
product of matrix rings over number fields containing $E$, in which 
case $G$, and thus (H1) and (H2), play no role.

The main motivation for this work has its origins in the following special case.
Let $L/K$ be a finite Galois extension of number fields with Galois group $G$
such that $E$ is a subfield of $K$ and put $d=[K:E]$. One can take $X$ to be 
$\mathcal{O}_{L}$, the ring of algebraic integers of $L$, 
and $\mathcal{A}$ to be the associated order
$$\mathcal{A}( E[G]; \mathcal{O}_{L}) 
:= \{ x \in E[G] \mid x(\mathcal{O}_{L}) \sseq \mathcal{O}_{L} \}.$$
The application of the algorithm to this special situation 
is implemented in Magma (\cite{magma}) under certain extra hypotheses
when $K=E=\Q$. The source code and input files are 
available from
\[ \texttt{http://www.mathematik.uni-kassel.de/$\sim$bley/pub.html}. \]
Other Galois modules to which the algorithm can be applied include
the $G$-stable ideals of $\mathcal{O}_{L}$ and, in certain cases, 
the torsion-free part of $\mathcal{O}_{L}^{\times}$.

The algorithm can be thought of as a non-abelian, higher rank 
generalization of the one given in \cite{bley-units};
though stated for the Galois module structure of units, this can be
adapted to general modules for $G$ abelian and $d=1$ with relatively few changes. 
It is also worth noting that under the same restrictions on $G$ and $d$,
the algorithm in \cite{bley-endres} computes the Picard group 
$\mathrm{Pic}(\mathcal{A})$ 
and solves the corresponding refined discrete logarithm problem,
thus computing a generator if it exists.

There is a considerable body of work related to the motivating special
case of the Galois module structure of rings of integers. We briefly mention just
a few of these results, using the notation above. The most progress has 
been made in the case that $L/K$ is at most tamely ramified.
In this setting, it is well-known that $\mathcal{A}=\mathcal{O}_{E}[G]$
and $\mathcal{O}_{L}$ is locally free over $\mathcal{O}_{E}[G]$ (see \cite{noether}).
The algorithm in \cite{bley-wilson} determines the class of 
$\mathcal{O}_{L}$ in the locally free class group $\Cl(\mathcal{O}_{E}[G])$,
and thus whether or not it is stably free (note that under hypothesis (H2)
all stably free $\mathcal{A}$-modules are in fact free). 
Important work of Fr\"ohlich and Taylor determines the class of $\OL$ in the 
locally free class group $\Cl(\Z[G])$ in terms of Artin root numbers of irreducible complex 
symplectic characters of $G$ (see \cite{frohlich-book}). Unfortunately, neither
of these approaches lead to any description of generators. However, explicit generators 
or algorithms to find them when $K=E=\Q$ and $G=A_{4}$, $D_{2p}$ ($p$ odd prime), $H_{8}$, $H_{12}$,  
or $H_{8} \times C_{2}$ are given in  \cite{cougnard-A4}, \cite{cougnard-D4}, \cite{martinet-D2p}, 
\cite{cougnard-queyrut-H12} and \cite{cougnard-H8-C2}.

When no assumption regarding the ramification of $L/K$ is made, the situation
is somewhat more difficult, not least because $\OL$ is not necessarily locally free over 
$\mathcal{A}$. Perhaps the most important result in this context is Leopoldt's
Theorem, which in the case that $K=E=\Q$ and $G$ is abelian shows that $\OL$ 
is always free over $\mathcal{A}$ and, in addition, explicitly constructs an element 
$\alpha \in \OL$ in terms of Gauss sums such that $\OL = \mathcal{A}\alpha$ 
(see \cite{leopoldt}; Lettl gives a simplified proof in \cite{lettl-global}). 
In the setting $K=E$ and $L/\Q$ abelian, progressively sharper
generalizations of Leopoldt's Theorem (with explicit generators) 
are given in \cite{cl}, \cite{bley}, \cite{byott-lettl} and \cite{me2}.

In future work, we hope to eliminate hypothesis (H2). Finding an algorithm
to explicitly compute Wedderburn decompositions and thereby eliminate
hypothesis (H1) is an independent problem in its own right, on which
some progress has been made by others.
A more detailed discussion of both hypotheses is given in Section \ref{alg-gens}.

\section{A Necessary and Sufficient Condition for Freeness}

Let $E$ be a number field and $G$ be a finite group.
Let $\mathcal{A}$ be any $\mathcal{O}_{E}$-order of full rank
in the group algebra $A := E[G]$, and let $\mathcal{M}$ be some
maximal $\mathcal{O}_{E}$-order in $A$ containing $\mathcal{A}$.
(In fact, the results of this section still hold when the group algebra $E[G]$ is
replaced by any finite-dimensional semisimple $E$-algebra.)
For any non-commutative ring $R$, we shall henceforth take ``$R$-module'' to
mean ``left $R$-module'', unless otherwise stated. 

If $\frp$ is a prime of $\mathcal{O}_{E}$ and $M$ is an $\mathcal{O}_{E}$-module,
we write $M_{\frp} := \mathcal{O}_{E, \frp} \tensor_{\mathcal{O}_{E}} M$ for the 
localization of $M$ at $\frp$. We say that $M$ is locally
free of rank $d$ if for every $\frp$, we have $M_{\frp}$ free 
over $\calA_\frp$ of rank $d$.
For an $\calA$-lattice $X$, i.e., an $\mathcal{A}$-module that is
finitely generated and torsion-free over $\mathcal{O}_{E}$, we set 
$\calM X := \calM \tensor_{\calA}X$ and usually identify $\calM X$ 
with the sublattice $\{ \lambda x \mid \lambda \in \calM, x \in X \}$ of the $E$-vector space $E \tensor_{\OE} X$. We define $\calM_\frp X_\frp$ in the same way.

The main results of this paper are consequences of the following proposition.

\begin{prop}\label{nec-suf-prop}
Let $X$ be an $\mathcal{A}$-lattice. Then $X$ is free of rank $d$ if and only if
\begin{enumerate}
\item $X$ is a locally free $\mathcal{A}$-lattice of rank $d$, and
\item there exist $\alpha_{1}, \ldots, \alpha_{d} \in X$ such that 
$\mathcal{M}X = \mathcal{M} \alpha_{1} \oplus \ldots \oplus \mathcal{M} \alpha_{d}$.
\end{enumerate}
Further, when this is the case, $X = \mathcal{A} \alpha_{1} \oplus \ldots \oplus 
\mathcal{A} \alpha_{d}$.
\end{prop}

\begin{proof}
If $X$ is a free $\mathcal{A}$-lattice of rank $d$ then (a) and (b) follow trivially. 

Suppose conversely that (a) and (b) hold and let 
$Y= \mathcal{A} \alpha_{1} \oplus \ldots \oplus \mathcal{A} \alpha_{d} \subseteq X$.
Both $X$ and $Y$ are locally free $\mathcal{A}$-lattices of rank $d$ and so 
for each non-zero prime $\frp$ of $\mathcal{O}_{E}$ there exists an isomorphism
$f_{\frp}: Y_{\frp} \lra X_{\frp}$ of $\mathcal{A}_{\frp}$-lattices which extends naturally
to an isomorphism $f_{\frp}: \mathcal{M}_{\frp}Y_{\frp} \lra \mathcal{M}_{\frp}X_{\frp}$ 
of $\mathcal{M}_{\frp}$-lattices. 
For each $\frp$ we have
$$ [ X_{\frp} : Y_{\frp} ]_{\mathcal{O}_{E, \frp}} =
 [ f_{\frp}(Y_{\frp}) : Y_{\frp} ]_{\mathcal{O}_{E, \frp}} =
 \mathrm{det}_{E}(f_{\frp})\mathcal{O}_{E,\frp},$$
where the two left-most terms are generalized module indices
(see \cite[II.4]{ft}). 
However, $\mathcal{M}Y = \mathcal{M}X$ and so each
$f_{\frp}:\mathcal{M}_{\frp}Y_{\frp} \lra \mathcal{M}_{\frp}X_{\frp} = \mathcal{M}_{\frp}Y_{\frp}$ is in fact an $\mathcal{M}_{\frp}$-automorphism and therefore also a
$\mathcal{O}_{E,\frp}$-automorphism. Hence 
$\mathrm{det}_{E}(f_{\frp}) \in \mathcal{O}_{E,\frp}^{\times}$ and so
$ [ X_{\frp} : Y_{\frp} ]_{\mathcal{O}_{E, \frp}} = \mathcal{O}_{E,\frp}$ for
each $\frp$. Together with the fact that $Y \sseq X$, this shows that 
$X=Y$.
\end{proof}

Let $R$ be a ring with identity and denote by $R^{\op}$ the opposite ring.
If $M$ is a free $R$-module of rank $d$, then a choice of basis for
$M$ induces an isomorphism
$\End_{R}(M) \cong \Mat_{d}(R)^{\op}$. Note that for any subring of a left Noetherian ring, 
there is no distinction between left and right multiplicative inverses or units 
(see \cite[Theorem 6.4]{reiner}).
Hence we have 
$\Aut_{R}(M) := \End_{R}(M)^{\times} 
\cong \GL_{d}(R)^{\op} := (\Mat_{d}(R)^{\op})^{\times}$ as groups. Since
$\GL_{d}(R)^{\op} = \GL_{d}(R)$ as sets, we shall henceforth drop the $^{\op}$
notation.

\begin{corollary}
Let $X$ be an $\mathcal{A}$-lattice. Then $X$ is  free  of rank $d$ if and only if
\begin{enumerate}
\item $X$ is a locally free $\mathcal{A}$-lattice of rank $d$,
\item there exist $\beta_{1}, \ldots, \beta_{d} \in \mathcal{M}X$ such that 
$\mathcal{M}X = \mathcal{M} \beta_{1} \oplus \ldots \oplus \mathcal{M} \beta_{d}$, and
\item there exists $\lambda \in \GL_{d}(\mathcal{M})^{}$ such that each $\alpha_{i} \in X$
where
$(\alpha_{1}, \ldots, \alpha_{d})^{\mathrm{T}} 
:= \lambda ( \beta_{1}, \ldots , \beta_{d} )^{\mathrm{T}}$.
\end{enumerate}
Further, when this is the case, 
$X = \mathcal{A}\alpha_{1} \oplus \ldots \oplus \mathcal{A}\alpha_{d}$.
\end{corollary}

Most of  the following notation is adopted from \cite{bley-boltje}. Denote the center of a ring $R$ by $Z(R)$. 
Set $C:=Z(A)$ and let $\mathcal{O}_{C}$ be the
integral closure of $\mathcal{O}_{E}$ in C. Let $e_{1}, \ldots, e_{r}$ be the primitive idempotents
of $C$ and set $A_{i}:=Ae_{i}$. Then
\begin{equation}
\label{weddecomp}
A = A_{1} \oplus \cdots \oplus A_{r}
\end{equation}
is a decomposition of $A$ into indecomposable ideals. Each $A_{i}$ is an $E$-algebra with identity 
element $e_{i}$. By Wedderburn's Theorem, the centers $E_{i} := Z(A_{i})$ are finite field extensions
of $E$ via $E \rightarrow E_{i}$, $\alpha \mapsto \alpha e_{i}$, and we have $E$-algebra isomorphisms
$A_{i} \cong \Mat_{n_{i}}(D_{i})$ where $D_{i}$ is a division ring with $Z(D_{i}) \cong E_{i}$. The
decomposition (\ref{weddecomp}) gives
\begin{equation}
C = E_{1} \oplus \cdots \oplus E_{r}, 
\quad \mathcal{O}_{C} = \mathcal{O}_{E_{1}} \oplus \cdots \oplus \mathcal{O}_{E_{r}},
\quad  \textrm{and} \quad \mathcal{M} =  \mathcal{M}_{1} \oplus \cdots \oplus \mathcal{M}_{r},
\end{equation}
where we have set $\mathcal{M}_{i}:=\mathcal{M}e_{i}$. This
in turn induces decompositions
\begin{eqnarray}
\Mat_d (\mathcal{M})^{} &=& \Mat_d (\mathcal{M}_{1})^{} \oplus \ldots \oplus 
\Mat_d(\mathcal{M}_{r})^{} 
\text{ and } \\
\GL_d (\mathcal{M})^{} &=& \GL_d (\mathcal{M}_{1})^{} \times \ldots \times 
\GL_d(\mathcal{M}_{r})^{}. 
\end{eqnarray}

For the rest of this section we suppose  $1 \leq i \leq r$ and $1 \leq j \leq d$.

\begin{corollary}\label{nec-suf-with-max-split}
Let $X$ be an $\mathcal{A}$-lattice. Then $X$ is free of rank $d$ if and only if
\begin{enumerate}
\item $X$ is a locally free $\mathcal{A}$-lattice of rank $d$,
\item for each $i$, there exist $\beta_{i,1}, \ldots, \beta_{i,d}$ such that
$\mathcal{M}_{i}X = \mathcal{M}_{i}\beta_{i,1} \oplus \ldots \oplus 
 \mathcal{M}_{i}\beta_{i,d}$, and
\item there exist $\lambda_{i} \in \GL_{d}(\mathcal{M}_{i})^{} $ such that each $\alpha_{j} \in X$, where 
$\alpha_{j}  := \sum_{i=1}^{r} \alpha_{i,j}$ \\ and
$(\alpha_{i,1}, \ldots , \alpha_{i,d})^{\mathrm{T}} := 
\lambda_{i} (\beta_{i,1}, \ldots , \beta_{i,d})^{\mathrm{T}}$.
\end{enumerate}
Further, when this is the case, 
$X = \mathcal{A}\alpha_{1} \oplus \ldots \oplus \mathcal{A}\alpha_{d}$.
\end{corollary}

Let $\mathfrak{f}$ be any full two-sided ideal of $\mathcal{M}$ contained in $\mathcal{A}$.
Then we have
$\mathfrak{f} \subseteq \mathcal{A} \subseteq \mathcal{M} \subseteq A$. 
Set $\overline{\mathcal{M}} := \mathcal{M}/\mathfrak{f}$ and 
$\overline{\mathcal{A}} := \mathcal{A}/\mathfrak{f}$ so that 
$\overline{\mathcal{A}} \subseteq \overline{\mathcal{M}}$ are finite rings, and denote the canonical
map $\mathcal{M} \rightarrow \overline{\mathcal{M}}$ by $m \mapsto \overline{m}$.
Note that we have decompositions
\begin{equation}
\mathfrak{f} = \mathfrak{f}_{1} \oplus \cdots \oplus \mathfrak{f}_{r} 
\quad  \textrm{and} \quad \overline{\mathcal{M}} =  
\overline{\mathcal{M}_{1}} \oplus \cdots \oplus \overline{\mathcal{M}_{r}},
\end{equation}
where each $\mathfrak{f}_{i}$ is a non-zero ideal of $\mathcal{M}_{i}$ and 
$\overline{\mathcal{M}_{i}} := \mathcal{M}_{i} / \mathfrak{f}_{i}$.

For each $i$, let $U_{i} \subset \GL_{d}(\mathcal{M}_{i})^{}$
denote a set of representatives of the image of the
natural projection $\GL_{d}(\mathcal{M}_{i})^{} \lra \GL_{d}(\overline{\mathcal{M}_{i}})^{}$.

\begin{corollary}\label{free-quotient-cor}
Let $X$ be an $\mathcal{A}$-lattice. Suppose that
\begin{enumerate}
\item $X$ is a locally free $\mathcal{A}$-lattice of rank $d$, and
\item for each $i$, there exist $\beta_{i,1}, \ldots, \beta_{i,d}$ such that
$\mathcal{M}_{i}X = \mathcal{M}_{i}\beta_{i,1} \oplus \ldots \oplus 
 \mathcal{M}_{i}\beta_{i,d}$.
\end{enumerate}
Then $X$ is free of rank $d$ over $\mathcal{A}$ if and only if
\begin{enumerate}
\setcounter{enumi}{2}
\item there exist $\lambda_{i} \in U_{i}$ such that each $\alpha_{j} \in X$, where 
$\alpha_{j}  := \sum_{i=1}^{r} \alpha_{i,j}$ \\ and
$(\alpha_{i,1}, \ldots , \alpha_{i,d})^{\mathrm{T}} := 
\lambda_{i} (\beta_{i,1}, \ldots , \beta_{i,d})^{\mathrm{T}}$.
\end{enumerate}
Further, when this is the case, 
$X = \mathcal{A}\alpha_{1} \oplus \ldots \oplus \mathcal{A}\alpha_{d}$.
\end{corollary}

\begin{proof}
If condition (c) holds, then the result follows immediately from Corollary \ref{nec-suf-with-max-split}.

Suppose conversely that $X$ is free of rank $d$ over $\mathcal{A}$. Then by
Corollary \ref{nec-suf-with-max-split} there exist
$\lambda_{i} \in \GL_{d}(\mathcal{M}_{i}^{}) $ such that each $\alpha_{j} \in X$
where the $\alpha_{j}$'s are defined as above. However, as $\mathfrak{f}$ is a 
two-sided ideal of $\mathcal{A}$, we have

\begin{eqnarray*} 
\bigoplus_{i=1}^{r} (\lambda_{i} + \Mat_{d}(\mathfrak{f}_{i})) (\beta_{i,1}, \ldots , \beta_{i,d})^{\mathrm{T}}
&\subseteq& \bigoplus_{i=1}^{r} \lambda_i  (\beta_{i,1}, \ldots , \beta_{i,d})^{\mathrm{T}} + 
\bigoplus_{i=1}^{r} \Mat_{d}(\mathfrak{f}_{i}) (\calM_i X)^d\\
&=& (\alpha_{1}, \ldots , \alpha_{d})^{\mathrm{T}} + 
\bigoplus_{i=1}^{r} \Mat_{d}(\mathfrak{f}_{i}) (\calM_i X)^d
\sseq X^{d}.
\end{eqnarray*}
Thus we can suppose without loss of generality that $\lambda_{i} \in U_{i}$ for each $i$.
\end{proof}

Let $L/K$ be a Galois extension of number fields with Galois group $G$ such that 
$E$ is a subfield of $K$. Let $d=[K:E]$ and write $\mathcal{O}_{L}$ for the ring of integers of $L$. 
One of the main applications of Corollary \ref{free-quotient-cor} is to determine whether
the ring of integers $\mathcal{O}_{L}$ is free of rank $d$ over the associated order
$\mathcal{A}=\mathcal{A}(E[G]; \mathcal{O}_{L}) := \{ x \in E[G] \mid x(\mathcal{O}_{L}) 
\sseq \mathcal{O}_{L}\}$.

In the case that $G$ is abelian and $E=K$, the maximal order $\mathcal{M}$ is unique
and everything can be made completely explicit in terms of the absolutely irreducible
characters of $G$. We refer the reader to \cite[Section 2.2]{bley-units}. 
The combination of Theorem 2.8 and Lemma 2.9 of loc. cit. is essentially equivalent to
Corollary \ref{free-quotient-cor} given here specialized to the abelian case.

We also remark that in the case that $E=K$ and $L/K$ is an at most 
tamely ramified Kummer extension with $G$ cyclic, results of Ichimura (see \cite[Theorem 2]{ichimura-kummer}) are, 
though not exactly the same, very similar to Corollary \ref{free-quotient-cor} when applied to this special situation.

\section{The Algorithm}\label{alg-gens}

Let $E$ be a number field and $G$ be a finite group.
Let $\mathcal{A}$ be any $\mathcal{O}_{E}$-order of full rank
in the group algebra $E[G]$ and let $X$ be an $\mathcal{A}$-lattice.
In this section, we give an algorithm based on Corollary \ref{free-quotient-cor} that 
either computes elements $\alpha_{1}, \ldots, \alpha_{d} \in \mathcal{O}_{L}$
such that 
$X = \mathcal{A} \alpha_{1} \oplus \ldots \oplus \mathcal{A} \alpha_{d}$,
or determines that no such elements exist. In other words, the algorithm 
determines whether $X$ is free over $\mathcal{A}$, and if so, computes
explicit generators. 

We require the hypotheses (H1) and (H2) formulated in the introduction, which we now recall and briefly remark upon. Note that the algorithm still works if the group algebra 
$E[G]$ is replaced by any finite product of matrix rings over number fields containing 
$E$, in which case $G$, and thus (H1) and (H2), play no role.

\begin{itemize}
\item [(H1)]  The Wedderburn decomposition $E[G]\cong \oplus_{\chi} M_{\chi}$,
where each $M_{\chi}$ is a matrix ring over a division ring, is explicitly computable.
\smallskip \\
If $G$ is abelian the Wedderburn decomposition can be explicitly computed from the character table. For $G$ non-abelian, many decompositions can be found in the literature or computed 
``by hand''. Note that this
problem is equivalent to explicitly finding all irreducible $E[G]$-modules up to isomorphism. 
An effective method that dates back to Schur to solve this important computational task
in the case where $G$ is soluble is likely to be implemented in Magma v2.14.
\smallskip
\item [(H2)] The Schur indices of all $E$-rational irreducible characters of $G$ are equal to 1,
i.e., each $M_{\chi}$ above is in fact a matrix ring over a number field. \smallskip \\
This holds, for example, whenever 
\begin{enumerate}
\item $G$ is abelian, dihedral or symmetric;
\item $G$ is a $p$-group where $p$ is an odd prime; or
\item $E$ contains a primitive $m$-th root of unity, where $m$ is the exponent of $G$.
\end{enumerate}
A full discussion of Schur indices is given in \cite[Chapter 10]{isaacs}. 
An algorithm of Nebe and Unger to compute the Schur index will be implemented in 
Magma v2.14 (a paper on this work is in preparation).
\end{itemize}

Before we sketch the individual steps of the algorithm, we briefly digress to describe
the presentation of our data. We always assume that $\calO_E[G]$-modules 
$X$ are given by an $\OE$-pseudo-basis as described, for example, 
in \cite[Definition~1.4.1]{cohen-advcomp}. 
To be more precise, we assume that $V := E \tensor_\OE X$ is
given by an $E$-basis $v_1, \ldots, v_m$ together with matrices $A(\sigma) \in \GL_m(E)$ 
for each $\sigma \in G$ describing the action of $G$,
\[
\left(
  \begin{array}{c}
    v_1\\ \vdots \\ v_m
  \end{array} \right)^\sigma = A(\sigma) \left( \begin{array}{c}
    v_1\\ \vdots \\ v_m
  \end{array} \right).
\]
Then $X = \fra_1 w_1 \oplus \ldots \fra_m w_m$, where each $\fra_i$ is a fractional ideal of 
$\OE$ and each $w_i \in V$. 
Similarly, $\calA = \frb_1 \lambda_1 \oplus \ldots \frb_n \lambda_n$ with fractional
$\OE$-ideals $\frb_i$ and $\lambda_i \in E[G]$.

\begin{algorithm}\label{thealg} 
Input: $\mathcal{A}$ and $X$ as above.
\renewcommand{\labelenumi}{(\arabic{enumi})}
\begin{enumerate}
\item Compute $d := \dim_{E} (E \otimes_{\mathcal{O}_{E}} X)/|G|$ and check that $d \in \N$.

\item Compute a maximal $\mathcal{O}_{E}$-order $\mathcal{M}$ in $E[G]$ containing $\mathcal{A}$.

\item Compute the central primitive idempotents $e_{i}$ and the components 
$\mathcal{M}_{i}:=\mathcal{M} e_{i}$.

\item Compute the conductor $\mathfrak{c}$ of $\mathcal{A}$ in $\mathcal{M}$ and the components
$\mathfrak{c}_{i} := \mathfrak{c} e_{i}$. \\ Then compute the ideals 
$\mathfrak{g}_{i} := \mathfrak{c}_{i} \cap \mathcal{O}_{E_{i}}$ and $\mathfrak{f}_{i} := \frg_{i} \mathcal{M}_{i}$ 
for each $i$.

\item For each $i$, compute $\beta_{i,1}, \ldots, \beta_{i,d}$ such that 
$\mathcal{M}_{i} X = 
\mathcal{M}_{i} \beta_{i,1} \oplus \ldots \oplus \mathcal{M}_{i} \beta_{i,d}$.

\item Check that $X$ is locally free of rank $d$ over $\mathcal{A}$.

\item For each $i$, compute a set of representatives 
$U_{i} \subset \GL_{d}(\mathcal{M}_{i})^{}$ 
of the image of the natural projection map 
$\GL_{d}(\mathcal{M}_{i})^{} \longrightarrow 
\GL_{d}(\overline{\mathcal{M}_{i}})^{}$,
where $\overline{\mathcal{M}_{i}} := \mathcal{M}_{i} / \mathfrak{f}_{i}$.

\item Find a tuple $(\lambda_{i}) \in \prod_{i=1}^{r} U_{i}$ such that 
that each $\alpha_{j} \in X$, where 
$\alpha_{j} := \sum_{i=1}^{r} \alpha_{i,j}$ \\ and
$(\alpha_{i,1}, \ldots , \alpha_{i,d})^{\mathrm{T}} := 
\lambda_{i} (\beta_{i,1}, \ldots , \beta_{i,d})^{\mathrm{T}}$.
For such a tuple, 
$X = \mathcal{A}\alpha_{1} \oplus \ldots \oplus \mathcal{A}\alpha_{d}$.

\end{enumerate}
\end{algorithm}

Before commenting on the individual steps, we remark that steps (1) to (4) 
can be done in full generality without assuming hypotheses (H1) or (H2).

\renewcommand{\labelenumi}{(\arabic{enumi})}
\begin{enumerate}

\item If we replace $E[G]$ by some finite product of matrix rings over number fields $A$, 
then we define $d := \dim_{E} (E \otimes_{\mathcal{O}_{E}} X)/ \dim_{E} (A)$. 

\item
An algorithm for computing $\mathcal{M}$ is described in
\cite[Kapitel 3 and 4]{friedrichs}.

\item
Each central primitive idempotent corresponds to an irreducible $E$-character $\chi_{i}$ 
and we have $e_{i} = \frac{n_{i}}{|G|} \sum_{g \in G} \chi(g^{-1})g$ with 
$n_{i}=\chi_{i}(1)$.

\item
In practice, we compute some multiple of the conductor. For example, one can use 
the method outlined in \cite[3.2 (f) and (g)]{bley-boltje}. Also see \cite[Remark 3.3]{bley-boltje}. 

\item
This step is described in Section \ref{max-over-numfields}.

\item 
Successful completion of step (5) shows that $\calM X$ is a free $\calM$-module
of rank $d$. Therefore $X$ is locally free of rank $d$ over $\mathcal{A}$
except possibly at the (finite number of) primes of $\mathcal{O}_{E}$ dividing the 
generalized module index $[\mathcal{M}:\mathcal{A}]_{\mathcal{O}_{E}}$.
An algorithm to compute local basis elements (and thus to check local freeness) 
at these primes is given in \cite[Section 4.2]{bley-wilson}. 
Note that in the motivating case $X=\mathcal{O}_{L}$ for some number field $L$ (see introduction), 
$\mathcal{M}X$ is always locally free over $\calM$  and so checking local freeness 
can be performed independently of step (5) and therefore without hypotheses (H1) or (H2). 
(To see this, note that $\calM X$ is projective over $\calM$ by \cite[Theorem~21.4]{reiner}, 
and $L$ is free over $K[G]$ and thus $E[G]$ by the Normal Basis Theorem.)

\item
This step is described in Section \ref{enunits}.

\item
The number of tests for this step can be greatly reduced 
by using a method analogous to the one outlined \cite[Section 2]{bley-units}.
We briefly describe this approach in Section \ref{improvement}.
However, even with this improvement, the enumeration is the most time-consuming 
part of the whole algorithm.
\end{enumerate}

\section{Computing Associated Orders}\label{assoc-orders}

Let $X$ be a finitely generated $\OEG$-module in the free $\EG$-space $V := E \tensor_\OE X$. 
In this section, we shall assume that an $\EG$-basis $v_1, \ldots, v_d$ of $V$ is known. 
The aim is to compute the order
\[
\calA(X) = \calA(E[G]; X):= \{ \lambda \in \EG \mid \lambda X \sseq X \}.
\]
We describe an algorithm which combines and contains all of the methods of \cite{bley-units}, 
\cite[Appendix]{burns-bley-appendix} and \cite[Lemma 3.1]{bley-endres}.

For further applications, such as the computation of conductors, we consider a more general problem and describe an algorithm to compute
\[
\calA(X, Y) = \calA(E[G]; X, Y) := \{ \lambda \in \EG \mid \lambda X \sseq Y \},
\]
where $Y \sseq V$ is another full $\OEG$-submodule.
Without loss of generality we may assume that $X, Y \sseq \EG^d$. 

We denote by 
$t : \EG \times \EG \lra E$ any symmetric, non-degenerate $E$-bilinear pairing. For computational
purposes we usually use the trace pairing which is characterized by
\begin{eqnarray*}
  t(g, h) =
  \begin{cases}
    1 & \text{if } gh = 1, \\ 0 & \text{otherwise,}
  \end{cases}\quad \text{ for } g, h \in G.
\end{eqnarray*}
We let $s : \EG^d \times \EG^d \lra E$ be the $d$-fold orthogonal sum of $t$. 
For any $\OEG$-module $M$ in $\EG^d$, respectively $\EG$, we identify the linear dual
$M^* := \mathrm{Hom}_\OE(M, \OE)$ with $\{ \lambda \in \EG^d \mid s(\lambda, M) \sseq \OE \}$, respectively $\{ \lambda \in \EG \mid t(\lambda, M) \sseq \OE \}$.
If $M$ is given by a pseudo-basis $(\mu_k, \frc_k)_{k}$, then $M^*$ is easy to compute. 
Indeed, if $\{\mu_k^*\}$ is the dual basis of $\{\mu_k\}$ with respect to $s$, respectively $t$, 
then $(\mu_k^*, \frc_k^{-1})_{k}$ is a pseudo-basis of $M^*$.
It is clear that the dual basis $\{\mu_k^*\}$ can be computed
by means of straightforward linear algebra. 

We now define an $\EG$-module homomorphism
\begin{eqnarray*}
( \cdot, \cdot )  : \EG^d \times EG^d &\lra& \EG, \\
(\mu, \nu) &\mapsto& \sum_{g \in G} s(g\mu, \nu)g^{-1}.
\end{eqnarray*}
This homomorphism satisfies 
\begin{equation}
  \label{eq300}
  t((\mu, \nu), \delta) = s(\nu, \delta \mu) =  s(\nu\delta, \mu)
\end{equation}
for $\mu, \nu \in \EG^d$ and $\delta \in \EG$.
 
\begin{lemma}
Let $V$ be a free $\EG$-space of rank $d$ and let $X, Y$ be two full $\OEG$-submodules of $V$.
Then $(X, Y^*) = \calA(X, Y)^*$.
\end{lemma}

\begin{proof}
Using (\ref{eq300}), this is essentially the same as the proof of
\cite[Lemma~4.2]{bley-burns-aa}.
\end{proof}

\begin{remark}
The main application is the following.
Let $L/K$ be a finite Galois extension of number fields with Galois group $G$
such that $E$ is a subfield of $K$ and put $d=[K:E]$. 
Let $I$ be an ambiguous (i.e. $G$-stable)  ideal of the ring of integers 
$\mathcal{O}_{L}$ and define the associated order to be 
$ \mathcal{A}(E[G];I) : =  \{ x \in E[G] \mid x(I) \sseq I \}.$
In \cite{girstmair-normal-basis}, an algorithm to compute a normal 
basis element for $L$ over $K$ (i.e. a generator for $L$ as a $K[G]$-module)\
is given, and from this it is easy to determine an $E[G]$-basis of $L$.
(It is also often easy to do this by trial and error.)
Hence we can apply the above method to compute the associated order and then, 
assuming hypotheses (H1) and (H2), find generators using Algorithm \ref{thealg}.
\end{remark}

\begin{remark}
The method of this section together with Algorithm \ref{thealg} can also be used to investigate 
the Galois module structure of units as in \cite{bley-units}. 
For a number field $L$, write $U_L$ for the units of $\mathcal{O}_{L}$ and $\mu(L)$ for the subgroup of roots of unity. Set $X := U_L / \mu(L)$ and write $A$ for the
semisimple algebra which acts naturally on $\Qu \tensor_\Z X$. 
The following cases can be considered:
\begin{itemize}
\item[(a)] $L/\Qu$ a totally real Galois extension, $A = \QG / \left( \sum_{g \in G} g \right)$;

\item[(b)] $L/\Qu$ a CM Galois extension with complex conjugation $\tau$, 
$A = \QG / \left( \tau - 1, \sum_{g \in G} g \right)$;

\item[(c)] $L/K$ a Galois extension of a quadratic imaginary field $K$, 
$A = \QG / \left( \sum_{g \in G} g \right)$.
\end{itemize}
Note that by \cite[Lemma~5.27]{wash} the module $\Qu \tensor_\Z X$ is free over $A$, so that $\calM X$ is
always locally free over $\calM$. Hence checking local freeness can be performed without the assumption
of hypotheses (H1) and (H2). 
\end{remark}

\begin{remark}
It is always possible to compute an $\EG$-basis $V = E \tensor_{\OE} X$ under hypotheses
(H1) and (H2) by using a weaker version of Proposition \ref{maxgen} in which the ring of integers
and its ideals are replaced by the appropriate number field.
\end{remark}

\section{Modules over Maximal Orders in Matrix Rings over Number Fields}\label{max-over-numfields}

Let $n \in \N$, let $F$ be a number field and let $\mathcal{O}=\mathcal{O}_{F}$ denote the ring of integers of $F$. 

\begin{prop}\label{maxform}
For each ideal $\mathfrak{a}$ of $\mathcal{O}$, let
$$ \Lambda_{\mathfrak{a},n} =
\left(\begin{array}{cccc}
\mathcal{O} & \cdots & \mathcal{O}  & \mathfrak{a}^{-1} \\
\vdots & \ddots & \vdots & \vdots  \\
\mathcal{O} & \cdots & \mathcal{O} & \mathfrak{a}^{-1}  \\
\mathfrak{a} & \cdots & \mathfrak{a}  & \mathcal{O}
\end{array}\right)$$
denote the ring of all $n \times n$ matrices $(x_{ij})$ where $x_{11}$ ranges over all elements of 
$\mathcal{O}$, \ldots, $x_{1n}$ ranges over all elements over $\mathfrak{a}^{-1}$, and so on. 
(For $n=1$, we take $\Lambda_{\mathfrak{a},n}=\mathcal{O}$.)
Then $\Lambda_{\mathfrak{a},n}$ is a maximal $\mathcal{O}$-order in $\Mat_{n}(F)$ and every maximal 
$\mathcal{O}$-order in $\Mat_{n}(F)$ is isomorphic to one of this form, for some ideal $\mathfrak{a}$ of $\mathcal{O}$.
\end{prop}

\begin{proof}
This is a special case of \cite[Corollary 27.6]{reiner}.
\end{proof}

Even though we can compute maximal orders (using \cite[Kapitel 3 and 4]{friedrichs}), 
we do not automatically get them in the above ``nice form''. 
We may assume that a maximal $\mathcal{O}$-order $\Lambda \subset \Mat_{n}(F)$
is given as an $\calO$-module by a $\calO$-pseudo basis. 
We briefly describe how to find an isomorphism that transforms $\Lambda$ into the ``nice form'' described in Proposition \ref{maxform}.

Let $Z \sseq F^n$ denote the $\calO$-module generated by the first column of $\Lambda$. 
Let
\[
Z = \calO z_1 \oplus \ldots \oplus \calO z_{n-1} \oplus \fra z_n, \quad z_i \in F^n,
\]
be the Steinitz form of $Z$ for some ideal $\mathfrak{a}$ of $\mathcal{O}$. 
(The Steinitz form of a torsion-free, finitely generated 
module over a Dedekind domain is the form given in \cite[Theorem 13(b)]{ft}.)

\begin{lemma}\label{general max order}
Let $S = \left( z_1, \ldots, z_n \right) \in \GL_n(F)$ be the matrix with columns $z_1, \ldots, z_n$. Then $\Lambda = S \Lambda_{\fra, n} S^{-1}$.
\end{lemma}

\begin{proof}
It is easy to see that $\Lambda = \{ \lambda \in \Mat_{n}(F) \mid \lambda Z \sseq Z \}$. 
With a slight abuse of
notation we may write $Z = (z_1, \ldots, z_n) (\calO, \ldots, \calO, \fra)^{\mathrm{T}} = 
S (\calO, \ldots, \calO, \fra)^{\mathrm{T}}$ and deduce
\begin{eqnarray*}
\lambda Z \sseq Z
&\iff& \lambda S (\calO, \ldots, \calO, \fra)^{\mathrm{T}} \sseq S (\calO, \ldots, \calO, \fra)^{\mathrm{T}} \\
&\iff& S S^{-1} \lambda S (\calO, \ldots, \calO, \fra)^{\mathrm{T}} \sseq S (\calO, \ldots, \calO, \fra)^{\mathrm{T}} \\
&\iff&  S^{-1} \lambda S (\calO, \ldots, \calO, \fra)^{\mathrm{T}} \sseq  (\calO, \ldots, \calO, \fra)^{\mathrm{T}} \\
&\iff& S^{-1} \lambda S \sseq \Lambda_{\fra,n}.
\end{eqnarray*}
\end{proof}

Replacing $\Lambda$ by $S^{-1} \Lambda S$ and a $\Lambda$-module $X$ by $S^{-1}X$ we may without loss of generality assume that our maximal order is in the above ``nice form''. 
We fix some maximal $\mathcal{O}$-order $\Lambda =  \Lambda_{\mathfrak{a},n}$ in 
$\Mat_{n}(F)$ for the rest of this section and now turn to the problem of determining whether a 
$\Lambda$-module $X$ is free of finite rank, and if so, whether generators can be computed.
Let $e_{kl}$ denote the matrix $(x_{ij}) \in \Lambda \subset \Mat_{n}(F)$ with $x_{ij}=0$ for $(i,j) \neq (k,l)$ and $x_{kl}=1$.

\begin{prop}\label{maxgen}
Let $X$ be a $\Lambda$-module. Then $X$ is free of rank $d$ over $\Lambda$, if and only if there exist 
$\omega_{1,1}, \ldots, \omega_{1,n}, \ldots , \omega_{d,1}, \ldots, \omega_{d,n} $ such that 
\[
e_{11}X = \mathcal{O}\omega_{1,1} \oplus \ldots \oplus  \mathcal{O} \omega_{1,n-1}
\oplus \mathfrak{a}^{-1} \omega_{1,n} \oplus \ldots \oplus
 \mathcal{O}\omega_{d,1} \oplus \ldots \oplus  \mathcal{O} \omega_{d,n-1}
\oplus \mathfrak{a}^{-1} \omega_{d,n}
\] 
Further, when this is the case, $X = \Lambda \omega_1 \oplus \ldots \oplus \Lambda \omega_d$
where $\omega_j := e_{11} \omega_{j,1} + \ldots + e_{n1} \omega_{j,n}$, $j =1, \ldots, d$.
\end{prop}

\begin{proof}
Suppose that $X$ is free of rank $d$ over $\Lambda$. 
Then $e_{11}$ ``cuts out the first row of each $\Lambda$'' in
$X \cong \oplus_{i=1}^{d} \Lambda$ and so $e_{11}X$ is of the desired form.

Now suppose conversely that there exist 
$\omega_{1,1}, \ldots, \omega_{1,n}, \ldots , \omega_{d,1}, \ldots, \omega_{d,n} $ such that 
\[
e_{11}X = \mathcal{O}\omega_{1,1} \oplus \ldots \oplus  \mathcal{O} \omega_{1,n-1}
\oplus \mathfrak{a}^{-1} \omega_{1,n} \oplus \ldots \oplus
 \mathcal{O}\omega_{d,1} \oplus \ldots \oplus  \mathcal{O} \omega_{d,n-1}
\oplus \mathfrak{a}^{-1} \omega_{d,n}
\] 
and define $\omega_j = e_{11} \omega_{j,1} + \ldots + e_{n,1} \omega_{j,n}$ for $j =1, \ldots, d$.

For $i \neq n$ and all $j$, we have $\omega_{j,i} \in e_{11}X \subset X$ and so $e_{i1}\omega_{j,i} \in X$. 
Furthermore, $e_{n1} \in \Lambda \mathfrak{a}^{-1}$ and 
$\omega_{j,n} \in \mathfrak{a}e_{11}X \subseteq \mathfrak{a} X$ for all $j$,
so $e_{n1} \omega_{j,n} \in X$. Therefore $\omega_{j} \in X$ for all $j$ and so 
$\Lambda \omega_{1} \oplus \ldots \oplus \Lambda \omega_{d} \subseteq X$.

Note that $X = e_{11}X \oplus \dots \oplus e_{nn} X$ since $e_{11} + \ldots + e_{nn}$ is the $n \times n$ 
identity matrix.
Furthermore, for all $j,k$ we have
$$ e_{1k} \omega_{j} = e_{1k}(e_{11} \omega_{j,1} + \ldots + e_{n1}\omega_{j,n}) 
= e_{1k}e_{k1}\omega_{j,k} = e_{11} \omega_{j,k} = \omega_{j,k}.$$
Therefore, since $\mathcal{O}e_{1k} \subseteq \Lambda$ for $k \neq n$ and $\mathfrak{a}^{-1}e_{1n}\subseteq \Lambda$, we have 
\begin{eqnarray*}
e_{11}X &=& \mathcal{O}\omega_{1,1} \oplus \ldots \oplus  \mathcal{O} \omega_{1,n-1}
\oplus \mathfrak{a}^{-1} \omega_{1,n} \oplus \ldots \oplus
 \mathcal{O}\omega_{d,1} \oplus \ldots \oplus  \mathcal{O} \omega_{d,n-1}
\oplus \mathfrak{a}^{-1} \omega_{d,n} \\
&=& \mathcal{O} e_{11} \omega_{1} \oplus \ldots \oplus
\mathcal{O} e_{1(n-1)}\omega_{1} \oplus \mathfrak{a}^{-1} e_{1n}\omega_{1} \oplus \ldots \oplus
\mathcal{O} e_{11} \omega_{d} \oplus \ldots \oplus
\mathcal{O} e_{1(n-1)}\omega_{d} \oplus \mathfrak{a}^{-1} e_{1n}\omega_{d} \\
&\subseteq& \Lambda \omega_{1} \oplus \ldots \oplus \Lambda \omega_{d}.
\end{eqnarray*}
Finally, observe that
\begin{eqnarray*}
e_{ii}X &=& e_{i1}e_{11}e_{1i}X \subseteq e_{i1}e_{11}X 
\subseteq e_{i1}( \Lambda \omega_{1} \oplus \ldots \oplus \Lambda \omega_{d}) \\ 
&\subseteq& \Lambda \omega_{1} \oplus \ldots \oplus \Lambda \omega_{d}
\quad \textrm{ for } i \neq n, \textrm{ and}\\
e_{nn}X &=& e_{n1}e_{11}e_{1n} X = (\mathfrak{a}e_{n1})e_{11}(\mathfrak{a}^{-1}e_{1n})X
\subseteq  (\mathfrak{a}e_{n1})e_{11} X \\
&\subseteq& (\mathfrak{a}e_{n1})
(\Lambda \omega_{1} \oplus \ldots \oplus \Lambda \omega_{d}) \subseteq 
\Lambda \omega_{1} \oplus \ldots \oplus \Lambda \omega_{d},
\end{eqnarray*}
so therefore 
$X = 
e_{11}X \oplus \dots \oplus e_{nn} X \subseteq \Lambda \omega_{1} \oplus \ldots \oplus \Lambda \omega_{d}$.
\end{proof}

\begin{corollary}\label{max-iso-cor}
Let $X$ be a $\Lambda$-module. Then $X$ is free of rank $d$ over $\Lambda$
if and only if $e_{11}X$ is of rank $dn$ and Steinitz class $[\mathfrak{a}^{-d}]$ as an 
$\mathcal{O}$-module.
\end{corollary}

We now give a description of Step (5) of Algorithm \ref{thealg}. Fix $i$, set 
$\Lambda=S^{-1}\mathcal{M}_{i} S$ and replace $X$ by $S^{-1}X$, 
where $S$ is as in Lemma \ref{general max order}. It is straightforward to see that it suffices to determine elements 
$\omega_{1,1}, \ldots, \omega_{d,n}$ satisfying the equation of Proposition \ref{maxgen}. 
First, compute a Steinitz form for $e_{11}X$, i.e. find $b_{j} \in e_{11}X$
and an ideal $\mathfrak{b}$ of $\mathcal{O}$ such that
$$ e_{11}X = \mathcal{O}b_{1} \oplus \ldots \oplus \mathcal{O} b_{dn-1} \oplus \mathfrak{b} b_{dn}$$
(one can use the Magma function \texttt{SteinitzForm}) and check that 
$[\mathfrak{b}] = [\mathfrak{a}^{-d}]$ in $\Cl(\mathcal{O})$.
Let $Y = \mathcal{O} b_{d(n-1)+1} \oplus \ldots \oplus \mathcal{O} b_{nd-1} \oplus \mathfrak{b} b_{dn}$
and compute $\mathfrak{a} Y$.  This is a free $\mathcal{O}$-module
of rank $d$ and so we can compute an $\mathcal{O}$-basis, $c_{1}, \ldots, c_{d}$, which
is also a ``$\mathfrak{a}^{-1}$ basis'' of $Y$. Now we can take $\omega_{j,n} = c_{j}$ for 
$j=1, \ldots, d$ and $\{ \omega_{j,k} \mid k \neq n \} = \{ b_{1}, \ldots, b_{d(n-1)} \}$.

\section{Enumerating Units}\label{enunits}

Let $d, n \in \N$, let $F$ be a number field and let $\mathcal{O}=\mathcal{O}_{F}$ denote the ring of integers of $F$. 
Let $\Lambda$ be some maximal $\mathcal{O}$-order of $\Mat_{n}(F)$. By Lemma 
\ref{general max order} we may assume that $\Lambda$ is of the ``nice form'' 
$\Lambda_{\fra, n}$. 
Let $\mathfrak{g}$ be some
non-zero ideal of $\mathcal{O}_{F}$ and let $\mathfrak{f} := \mathfrak{g}\Lambda$. 
Throughout this section, we identify $\Mat_{d}(\Lambda)$ with a subring of $\Mat_{dn}(F)$ in the obvious way.
We wish to compute a set of
representatives $U \subset \GL_{d}(\Lambda)$ of the image of the natural projection map 
$\pi : \GL_{d}(\Lambda) \longrightarrow \GL_{d}(\overline{\Lambda})$ 
where $\overline{\Lambda} = \Lambda / \mathfrak{f}$.

\begin{definition}
Let $i,j \in \{1, \ldots, nd \}$ with $i \ne j$ and let
$$ 
x \in 
\begin{cases} 
\calO / \frg, & \text{ if } i,j \nmid n \text{ or } i,j \mid n, \\
\fra^{-1} / \frg\fra^{-1}, & \text{ if } i \nmid n \text{ and } j \mid n, \\
\fra / \frg\fra, & \text{ if } j \nmid n \text{ and }  i \mid n.
\end{cases}
$$
Then the elementary matrix $\E_{ij}(x)$ is the  matrix in $\GL_d(\overline{\Lambda})$ that 
has $1$ in every diagonal entry, has $x$ in the $(i,j)$-entry and is zero elsewhere.
Let $\E(\overline{\Lambda})$ denote the subgroup of $\GL_{d}(\overline\Lambda)$ 
generated by all elementary matrices and define $\E(\Lambda)$ analogously.
Note $\pi(\E(\Lambda)) = \E(\overline\Lambda)$.
\end{definition}

\begin{prop}\label{gen-ev}
Let $\varepsilon:\mathcal{O}_{F}^{\times} \rightarrow (\mathcal{O}_{F}/\mathfrak{g})^{\times}$
be the natural projection map and let $V$ be the subgroup of matrices
$(x_{ij}) \in \GL_{d}(\overline{\Lambda})$
with $x_{11} \in \varepsilon(\mathcal{O}_{F}^{\times})$, $x_{ii}=1$ for $i \neq 1$ 
and $x_{ij}=0$ for $i \neq j$. Then $\pi(\GL_{d}(\Lambda))$ is generated by 
$\E(\overline{\Lambda})$ and $V$. 
\end{prop}

\begin{proof}
Let $R := \calO / \frg$ and consider the $R$-modules 
$X := \oplus_{i=1}^{d} R^n$ and 
$Y := \oplus_{i=1}^{d} (R^{n-1} \oplus \fra/\frg\fra)$.
We choose $\xi \in F^\times$ and an integral ideal $\frb$ such that
\[
\fra = \xi\frb, \quad \frb + \frg = \calO.
\]
Let $b \in \frb$ and $y \in \frg$ such that $b + y = 1$. Then we have an isomorphism 
$\calO / \frg \lra \fra / \fra\frg$ of $R$-modules defined by $z + \frg \mapsto zb\xi + \fra\frg$. The inverse
is given by $z+\fra\frg \mapsto \xi^{-1}z + \frg$.

This induces an isomorphism $\varphi \colon X \ra Y$, and as a consequence we obtain an isomorphism
\begin{eqnarray*}
\psi \colon \GL_{nd}(R) &\lra& \GL_{d}(\overline{\Lambda}), \\
\overline{A} = (\overline{A}_{ij})_{1 \leq i,j \leq d} &\mapsto& (\overline{\Phi_{2} A_{ij} \Phi_{1}})_{1 \leq i,j \leq d} 
\end{eqnarray*}
where $A_{ij} \in \Mat_{n}(\calO)$,
$$
\Phi_1 = \left( \begin{array}{cccc} 1 & & & \\ & \ddots && \\ &&1& \\ &&&\xi^{-1} \end{array} \right)
\text{  and } 
\Phi_2 = \left( \begin{array}{cccc} 1 & & & \\ & \ddots && \\ &&1& \\ &&&b\xi \end{array} \right).
$$

One easily verifies that $\psi(\E_{nd}(R)) = \E(\overline\Lambda)$ where $\E_{nd}(R)$ denotes the
group generated by elementary matrices of $\Mat_{nd}(R)$. From \cite[Corollary (9.3), p.267]{bass} we
deduce $\SL_{nd}(R) = \E_{nd}(R)$. Hence we have a commutative diagram with exact rows 

$$
\xymatrix@1@!0@=48pt { 
1  \ar[r] & \SL_{nd}(R) \ar[rr] \ar[d]^\simeq_\psi & & \GL_{nd}(R) \ar[rr]^{\det}  
\ar[d]_{\psi}^\simeq & & R^{\times} \ar[r] \ar[d]_{=} &  1 \\
1 \ar[r] & \E(\overline\Lambda) \ar[rr] & & \GL_{d}(\overline{\Lambda}) \ar[rr]^{\det'}  
& & R^\times \ar[r] & 1
}
$$
where $ \det' := \det \circ \psi^{-1}$. The diagram
$$
\xymatrix@1@!0@=48pt { 
& & & \GL_{d}(\Lambda) \ar[rr]^{\det}  
\ar[d]^\pi & &  \mathcal{O}^\times \ar[r] \ar[d]^\varepsilon & 1 \\
1 \ar[r] & \E(\overline\Lambda) \ar[rr] & & \GL_{d}(\overline{\Lambda}) \ar[rr]^{\det'}  
& &  R^\times \ar[r] & 1
}
$$
also has exact rows and a straightforward computation shows that it commutes.
This immediately implies the assertions of the proposition.
\end{proof}

We now give a description of Step (7) of Algorithm \ref{thealg}. 
Fix $i$,  and set
$n=n_{i}$, $\Lambda=S^{-1}\mathcal{M}_{i}S$ with $S$ as in Lemma \ref{general max order}, 
$\mathfrak{g}=\mathfrak{g}_{i}, F=E_{i}$ and $U=U_{i}$.
Using, for example, \cite[Algorithm~6.5.8]{cohen-comp}, compute a generating set $\{ a_{1}, \dots, a_{s}\}$
for $\mathcal{O}_{F}^{\times}$. Then $\{ \varepsilon(a_{1}), \ldots, \varepsilon(a_{s}) \}$
is a generating set for $\varepsilon(\mathcal{O}_{F}^{\times})$ and using the obvious
isomorphism we have a generating set for $V$. The group $\E(\Lambda)$ is generated
by the elementary matrices $\E_{ij}(b_{ijk})$ for $i,j \in \{1, \ldots, n \}$, $i \neq j$ where for
fixed $i,j$, $\{ b_{ijk} \}$ is a $\Z$-spanning set for $\mathcal{O}/\mathfrak{g}$,
$\mathfrak{a}/\mathfrak{ga}$ or $\mathfrak{a}^{-1}/\mathfrak{ga}^{-1}$, as appropriate.
Such spanning sets can be computed using Hermite Normal Form techniques described, for
example, in \cite[Chapter~2.4]{cohen-comp}. By Proposition \ref{gen-ev}, we now have an 
explicit generating set for $\pi(\GL_{d}(\Lambda))$, and so it is straightforward to compute the 
desired set of representatives $U=U_{i}$.

\section{Reducing the Number of Final Tests}\label{improvement}

The final number of tests in step (8) of Algorithm \ref{thealg} can be enormous. For example, if $G \simeq S_4$ (the symmetric group with $24$ elements) and $\calA = \ZG$, then a computation shows that there are approximately $4.4 \times 10^{18}$ tuples 
$(\lambda_i) \in \prod_{i=1}^r U_i$, which need to be tested. In this section, we describe 
an ad hoc method analogous to the one outlined in \cite[Section 2]{bley-units} to reduce the
number of tests required.

However, even with this improvement, the number of tests which need to be performed is still
very large. Despite this, somewhat surprisingly, we can find generating elements in 
many $S_4$-examples. It would be interesting to have an explanation, possibly probabilistic or heuristic in nature, for this phenomenon.

The improvement is based on the following simple observation. Let
\begin{eqnarray*}
  \calM X &=& \fra_1 v_1 \oplus \ldots \oplus \fra_m v_m, \\
  X &=& \frb_1 w_1 \oplus \ldots \oplus \frb_m w_m,
\end{eqnarray*}
be $\OE$-pseudo-basis representations of $\calM X$ and $X$. Let $A \in \GL_m(E)$ be the transformation matrix such that
\[
\left(
  \begin{array}{c}
    w_1\\\vdots\\w_m   \end{array}\right) = A \left(
  \begin{array}{c}
    v_1\\\vdots\\v_m   \end{array}\right).
\]

We now apply the Hermite Normal Form algorithm in Dedekind domains (see
\cite[Algorithm 1.4.7]{cohen-advcomp}) to the matrix $A$ and the list of ideals 
$(\frb_1, \ldots, \frb_m)$, though we reduce rows rather than columns.
We obtain a matrix $U \in \GL_m(E)$ and a list of ideals $(\frc_1, \ldots, \frc_m)$ such that the matrix $H = UA$ is upper triangular with $1$ on each diagonal entry. Moreover,
\[
\frc_1 h_1 \oplus \ldots \oplus \frc_m h_m = \frb_1 a_1 \oplus \ldots \oplus \frb_m a_m,
\]
where $h_1, \ldots, h_m$ denote the rows of $H$ and $a_1, \ldots, a_m$ denote the rows of $A$.
This immediately implies that $U (w_1, \ldots, w_m)^{\mathrm{T}}$ together with the list of ideals 
$(\frc_1, \ldots, \frc_m)$ is also a pseudo-basis for $X$. 

Now suppose that the vector $(x_1, \ldots, x_m) \in E^m$ defines an element 
$x = \sum_{i=1}^m x_i v_i \in \calM X$. Then we have
\begin{equation}
  \label{eq301}
  x \in X \iff (x_1, \ldots, x_m) H^{-1} \in (\frc_1, \ldots, \frc_m).
\end{equation}
Since $H^{-1}$ is upper triangular, this leads to a much more efficient enumeration. In addition, in many cases the coefficients $x_1, \ldots, x_m$ can be easily computed by a clever choice
of basis $v_1, \ldots, v_m$. To illustrate this, we conclude this section with a brief discussion of the case where $G \simeq S_n$, $E = \Q$ and $X \sseq E[G]$ is locally free of rank $1$.

Let 
\[
\Phi : \QG \lra \bigoplus_{i=1}^r \Mat_{n_i} (\Q)
\] 
be the explicitly computable isomorphism that gives the Wedderburn decomposition of $\QG$.
Let $\calM \sseq \QG$ be the maximal order such that $\Phi(\calM)=\oplus_{i=1}^r \Mat_{n_i}(\Z)$. 
For reasons of efficiency, we choose to work with matrices and henceforth consider
$\mathcal{M}X$ as a module over $\oplus_{i=1}^r \Mat_{n_i}(\Z)$ via the isomorphism $\Phi^{-1}$.
Let $B = (B_1, \ldots, B_r)$ denote a $\oplus_{i=1}^r \Mat_{n_i}(\Z)$-basis of $\mathcal{M}X$.
Let $e_{i, kl} = (\ldots, e_{kl}, \ldots )\in \oplus_{j=1}^r\Mat_{n_j}(\Z)$, 
$i=1, \ldots, r$, $1 \le k,l \le n_i$, denote the tuple of matrices with the matrix 
$e_{kl}$ in the $i$-th position and the zero matrix everywhere else. 
Then the set $\{ e_{i, kl} B \}$ forms a $\Z$-basis of $\calM X$.

Now let $(\lambda_1, \ldots, \lambda_r) \in \prod_{i} U_i$. Then the coefficients $(x_{i,kl})_{i,k,l}$ of $(\lambda_1 B_1, \ldots, \lambda_r B_r)$ with respect 
to the basis $\{ e_{i, kl}B \}$ are given by the coefficients of the matrices $\lambda_i$ because
\[
\lambda_i B_i = \left( \sum_{1 \le k,l \le n_i} \lambda_{i,kl}e_{kl} \right) B_i = 
 \sum_{1 \le k,l \le n_i} \lambda_{i,kl} \left( e_{kl} B_i \right).
\]

\section{Implementation and Computational Results}

In this section, we describe the cases for which Algorithm \ref{thealg} has been 
implemented in Magma (\cite{magma}).
The source code and input files are available from
\[ \texttt{http://www.mathematik.uni-kassel.de/$\sim$bley/pub.html}. \]

Let $L/K$ be a finite Galois extension of number fields with Galois group $G$
such that $E$ is a subfield of $K$ and put $d=[K:E]$. As discussed in Section
\ref{assoc-orders}, Algorithm \ref{thealg} can be applied in this situation with 
$X=\mathcal{O}_{L}$ and $\mathcal{A}  = \mathcal{A}( E[G]; \mathcal{O}_{L})$. 
However, for the sake of simplicity, all aspects of the implementation in Magma 
are restricted to the case $K=E=\Q$.

Let $\mathcal{A}_{L/\Q} = \mathcal{A}(\Q[G]; \OL)$.
We have the following:
\renewcommand{\labelenumi}{(\alph{enumi})}
\begin{enumerate}

\item For any finite Galois extension $L/\Q$, we can compute the associated order 
$\mathcal{A}_{L/\Q}$ 
and check that $\mathcal{O}_{L}$ is locally free over $\mathcal{A}_{L/\Q}$, provided that Magma 
can compute the ring of integers $\mathcal{O}_{L}$ and that the Magma function
$\texttt{AutomorphismGroup(L)}$ works. Of course, this can be improved if theoretical information 
for either the ring of integers or the Galois group is available.

\item For $G=A_{4},S_{4}, D_n$ or $G$ abelian, we can explicitly compute the Wedderburn decomposition of 
$\Q[G]$ so that hypothesis (H1) is satisfied (here $D_{n}$ is the dihedral group of order $2n$).

\item We can compute generators $\alpha_{i}$ such that 
$\mathcal{M}_{i} \mathcal{O}_{L} = \mathcal{M}_{i} \alpha_{i}$ whenever
$G=A_{4},S_{4}$, $D_{n}$ or $G$ abelian. This works very well for small $n$ and small abelian groups. 
For example, we successfully ran many 
experiments with dihedral groups $D_n$ and $n \le 10$.
Note however, that our implementation requires that all the fields $E_i$ have class number $1$. 

\item We can compute a generator $\alpha$
such that $\mathcal{O}_{L} = \mathcal{A}_{L/\Q} \alpha$ whenever $G=A_{4}$, $D_{n}$ 
with $n$ small or $G$ a small abelian group. For dihedral groups ``small'' means something
like $n \le 10$, for abelian groups experiments show that we can easily deal with groups of order
$\le 20$. For $S_4$-extensions the number of checks required in the final enumeration 
is simply too large to be done in a naive way. Indeed, there are five Wedderburn components of 
$\Q[S_4]$ and if $L/\Qu$ is tame, then the numbers of elements in the sets $U_i$ 
(notation as in Proposition \ref{free-quotient-cor}) are
\[
|U_1| = 2, \quad |U_2| = 2, \quad |U_3| = 2304, \quad |U_4| =  22020096, \quad 
|U_5| =  22020096.
\]
Note that these numbers are smaller if $L/\Q$ is wildly ramified.

The authors implemented the reduction method outlined in Section \ref{improvement}
and, to their surprise, were able to compute a generator in all of the examples tested 
for which $\OL$ is locally free over its associated order. (The fact that a generator exists in 
this situation is not surprising - see discussion below). As already mentioned, it would be very interesting to have some explanation for this unexpected phenomenon. Furthermore,
$\OL$ failed to be locally free over its associated order in all the computed $S_{4}$ examples
for which $L/\Q$ is wildly ramified at both $2$ and $3$, though no examples of this were found when only one prime is wildly ramified.

\item The algorithm as implemented in Magma is not deterministic, i.e., the program
will produce different generators for the same extension when run at different times.
The relevant steps, where different choices may finally lead to different generators, are the
choice of a normal basis element for $L/\Qu$ and the order of the final enumeration. 

\end{enumerate}

The authors computed generators for more than $140$ extensions $L/\Q$
with Galois group $A_{4}$ taken from the tables of \cite{klueners-malle}.
This might lead one to speculate, for example, that every such extension
has the property that $\mathcal{O}_{L}$ is free over $\mathcal{A}_{L/\Q}$. 
In principle, one can prove or disprove this assertion in the following way. 

It is well-known that the locally free class group
$\Cl(\Z[A_{4}])$ is trivial (see \cite{cougnard-A4}, for example), and from this it is straightforward to show
that $\Cl(\mathcal{A})$ is also trivial for any order $\mathcal{A}$ with 
$\Z[A_{4}] \sseq \mathcal{A} \sseq \Q[A_{4}]$. Since $\Q[A_{4}]$ satisfies the
Eichler condition relative to $\Z$ (see \cite[Definitions 34.3 and 38.1]{reiner}), a result of Jacobinski shows that if 
an $\mathcal{A}$-module has trivial class in $\Cl(\mathcal{A})$, then it is in fact free over $\mathcal{A}$ 
(see \cite[Theorem 38.2]{reiner}, for example). Hence we are reduced to establishing whether $\mathcal{O}_{L}$ 
is locally free over $\mathcal{A}_{L/\Q}$ for every $A_{4}$-extension $L/\Q$. The authors thank the referee for the following observation: as the ramification filtrations for primes above $2$ or $3$ are very restricted, it seems plausible that one could in fact determine all possible associated orders and local Galois module structures by hand (as far as the authors and the referee are aware, no-one has actually done this). One might also be able to carry out a complete analysis for $S_{4}$-extensions of $\Q$, where again it is known that the locally free class group is trivial (see \cite{reiner-ullom-74}) and $\Q[S_{4}]$ satisfies the
Eichler condition relative to $\Z$. (However, unlike the $A_{4}$ case, there are many 
known examples of $S_{4}$-extensions for which local freeness fails - see (d) above.)

In fact, for any number field $K$ and any finite group $G$, it is possible to check in a finite amount of time whether 
every extension $L/K$ with Galois group $G$ has the property that $\mathcal{O}_{L}$ is locally free over 
$\mathcal{A}_{L/K} := \mathcal{A}(K[G] ; L)$. By Noether's Theorem (see \cite{noether}), we have local freeness at all primes of $\mathcal{O}_{K}$ that are at most tamely ramified in $L/K$. Therefore, it suffices to check all extensions of $p$-adic fields
$K_{\mathfrak{p}}$ with Galois group $H$, where $\mathfrak{p}$ ranges over all primes of $\mathcal{O}_{K}$
dividing the order of $G$ and $H$ ranges over all solvable subgroups of $G$. Enumerating all such extensions
is possible using the algorithm of \cite{pauli-roblot} (this gives generating polynomials for all extensions of a $p$-adic field $K_{\mathfrak{p}}$ of given degree and discriminant), and local freeness can be checked using
the method outlined in \cite[Section 4.2]{bley-wilson}.

\section{Acknowledgments}\label{acknowledgments}

The authors are grateful to the Deutscher Akademischer Austausch Dienst (German Academic Exchange Service) for a grant allowing the second named author to visit Cornelius Greither at Universit\"at der Bunderswehr M\"unchen 
for the 2006-07 academic year, thus making this collaboration possible. The authors also 
wish to thank John Cannon for useful correspondence regarding features to appear in the
next version of Magma, and the referee for several helpful comments and suggestions.

\bibliography{ComputeGeneratorsBib}{}
\bibliographystyle{amsalpha}

\end{document}